\documentclass[a4paper,12pt]{article}
\usepackage{amssymb,amsthm,latexsym}
\newcommand{\uh}{\underline{h}}
\newcommand{\Med}[1]{\left\langle #1 \right\rangle}
\newcommand{\med}[1]{\langle #1 \rangle}
\newcommand{\vol}{\Lambda_N}

\newcommand{\EE}{{\mathbb E}}
\newtheorem{theorem}{Theorem}

\newtheorem{proposition}{Proposition}

\begin{document}

\title{The Ising-Sherrington-Kirpatrick model 
in a magnetic field at high temperature}
\author{Francis~Comets% 
\thanks{{  University Paris 7 - Denis Diderot}.
Partially supported by  CNRS,
UMR 7599
``Probabilit\'es et Mod\`{e}les Al\'eatoires." }
\and Francesco Guerra\thanks{ Universit\`a di Roma ``La Sapienza'' and INFN, Sezione di Roma. 
Partially supported by MIUR and INFN
}
\and
Fabio Lucio Toninelli
\thanks{Laboratoire de Physique, Ecole Normale Sup\'erieure de Lyon. On leave at Institut f\"ur Mathematik,
Universit\"at Z\"urich
}}

%%%%%%%%%%%%%%%%%%%%%%%%%%%%%%%%%%%%%%%%  date
\date{December 2004} 

\maketitle

%%%%%%%%%%%%%%%%%%%%%%%%%%%%%%%%%%%%%%%%%%%%
\begin{abstract}
We study a spin system on a large box with both Ising interaction
and Sherrington-Kirpatrick couplings, in the presence of an external field.
Our results are: (i) existence of the pressure in the limit of 
an infinite box.
When both  Ising and  Sherrington-Kirpatrick temperatures are 
high enough, we prove that:
(ii)  the value of the pressure is given by a suitable  replica 
symmetric solution, and
(iii) the fluctuations of the pressure 
are of order of the inverse of the square of the volume with a 
normal distribution in the limit.
In this regime,  the pressure can be expressed in terms of 
random field Ising models. 
\end{abstract}
\noindent
{\small {\bf Key Words:} Ising model, Sherrington-Kirpatrick model,
spin-glass, thermodynamic limit, pressure, quadratic coupling\\
{\bf AMS (1991) subject classifications:} 
82B44; 60K35, 82D30
\\
{\bf Short title:} Ising-Sherrington-Kirpatrick model }
%%%%%%%%%%%%%%%%%%%%%%%%%%%%%%%%%%%%%%%%%%%%

\pagestyle{myheadings}
\markboth{}{Ising-Sherrington-Kirpatrick model}

\section{Introduction}

We consider a  $d$-dimensional Ising model in a magnetic field $h$,
perturbed by a mean field interaction of spin-glass type. 
The Hamiltonian contains two parameters, $\beta$ and $\kappa$, which play the role of two inverse temperatures.
When $\beta=0$ the model reduces to the Ising model at temperature $1/\kappa$, while for $\kappa=0$ one recovers
the Sherrington-Kirkpatrick (SK) model at temperature $1/\beta$.
The understanding of the SK model has recently witnessed great progress (see, e.g., \cite{GueTon02pres}
\cite{broken} 
\cite{taladef}).
The main interest in the analysis of this model is the possibility of investigating the robustness
of the phenomena typical of mean field spin glass models, in the presence of additional interactions
of non-mean-field character.

The model has been previously considered in \cite{ComGiaLeb} under the additional assumptions that $h=0$ and 
that the Ising model is ferromagnetic. Under these conditions it was proven 
that, if $\beta$ and 
$\kappa$ are small enough, the infinite volume pressure is given by the sum of the Ising pressure and of the 
SK one, at the respective temperatures.
Moreover, the disorder fluctuations of the pressure were found to be of order $1/V$, $V$ being the volume of the
system, and to satisfy a central limit theorem. 

In the present paper, the Ising interaction decays exponentially fast with distance,
but is not necessarily ferromagnetic. It turns out that
the presence of the magnetic field changes qualitatively the picture with respect to \cite{ComGiaLeb}. 
Indeed we find, still
for  $\kappa$ and $\beta$ small enough, that the limit pressure is given in terms of the pressure of an 
Ising model with random external field, the strength of the randomness being related to the typical value of the
overlap between two replicas of the system. 
It is remarkable, though natural, that the random field Ising model,
which has its own interest \cite{AizWeh90}, plays such an important
role in our model. Note also that the pressure of our 
model can be computed only via the thermodynamic limit of another 
disordered system. This is contrast with the case 
$h=0$ we mentioned above, and also of course 
with the case $\kappa =0$ of the standard
SK model. 
The second difference is that the fluctuations of the pressure in presence of $h$ satisfy a central limit
theorem on the scale $1/\sqrt V$ rather than $1/V$. The same phenomenon is known to happen in the 
case of the usual SK model, see for instance \cite{alr} \cite{cn} \cite{GueTon02tcl}.

In the region of thermodynamic parameters we consider, the system is in a ``replica-symmetric'' (RS) phase, the 
overlap between two independent replicas being a non-random value in the thermodynamic limit.
%, and correlations between local functions
%decay exponentially with the distance \footnote{(******in principle, we should prove the statement about 
%correlations. 
%I think it is easy, since we know it for 
%$t=0$ and this property should be conserved (in mean square sense) for t up to 1. However, adding the 
%proof could make the paper pedantic. maybe we could add a remark (without proof) about this fact after theorem 2,
%if we are sure it is correct. to be decided********* F: How would you 
%prove that property is conserved ? ******)}.
Our methods fail beyond some values $\beta_0(h)$ and $\kappa_0(h)$, which we believe
to be an artifact of our approach, rather than representing the true boundary of the RS region. 
The same inconvenient has been previously encountered in the analysis of the SK model \cite{T02}, \cite{GueTon02quad}.
In principle one could
improve these values by employing the replica symmetry breaking scheme which for instance 
enabled M. Talagrand \cite{T} 
to control the whole replica-symmetric region of the SK model, and later 
the entire phase space \cite{taladef}. However, in the present case (as well as in
\cite{ComGiaLeb}) one of the reasons
why we do not reach the true critical line is due to an incomplete control of the underlying random field Ising model,
and this problem would not be fixed by the methods of \cite{T}.
For this reason, we prefer to use a generalization of the technically simpler ``quadratic
replica coupling'' technique introduced in \cite{GueTon02quad}.

It would be of course an interesting challenge to go beyond the present approach, and to deal with 
lower-temperature situations where the Ising-SK system possibly shows a RSB-like behavior.

\section{Description of the model and results}

The model we consider is defined on the  $d$-dimensional hypercubic box 
$\Lambda_N=\{-N,\cdots,N\}^d$ and
its partition function is:
\begin{eqnarray}
\label{H}
Z_N(\kappa,\beta,h;J)=\sum_{\sigma\in\{-1,+1\}^{\Lambda_N}}\exp\left(-
\kappa H^I_N(\sigma)-\beta H^{SK}_N(\sigma;J)+h \sum_{i\in \Lambda_N}\sigma_i\right).
\end{eqnarray}
The Hamiltonian of the Sherrington-Kirkpatrick (SK) model is defined as
\begin{eqnarray*}
H^{SK}_N(\sigma;J)=-\frac1{\sqrt{2|\vol|}}\sum_{i,j\in\vol}J_{ij}\sigma_i\sigma_j
\end{eqnarray*}
and the couplings $J_{ij}$ are i.i.d. Gaussian random variables ${\mathcal N}(0,1)$.
On the other hand, the Hamiltonian of the Ising model is 
\begin{eqnarray*}
H^I_N(\sigma)=-\frac12 \sum_{i,j\in\vol} K(i-j)\sigma_i\sigma_j,
\end{eqnarray*}
where we assume that the interaction decays exponentially, i.e., 
\begin{equation}
\label{expdec}
|K(i)|\le C_1 e^{-C_2 |i|}
\end{equation}
for some $C_1,C_2>0$.
We do not require the interaction to be ferromagnetic.
The finite volume (disorder-dependent) pressure is defined as usual as
\begin{eqnarray*}
p_N(\kappa,\beta,h;J)=\frac1{|\vol|}\log Z_N(\kappa,\beta,h;J).
\end{eqnarray*}

Later, we will need to consider the (Gaussian) Random Field Ising Model, defined by
the partition function
\begin{eqnarray*}
Z^{\rm  RFIM}_N(\kappa,h,\gamma;J)=\sum_{\sigma\in\{-1,+1\}^{\Lambda_N}}\exp\left(-
\kappa H^I_N(\sigma)+\sum_{i\in \Lambda_N}\sigma_i(h+{\gamma} J_i)\right),
\end{eqnarray*}
$J_i$ being i.i.d. standard Gaussian variables, and being also independent
from the $J_{ij}$'s in the sequel. The existence of
the infinite volume pressure of the RFIM,
\begin{eqnarray*}
p^{\rm  RFIM}(\kappa,h,\gamma)&=&\lim_{N\to\infty} \EE\, p^{\rm  RFIM}_N(\kappa,h,\gamma;J)\\
&=&a.s.-\lim_{N\to\infty} p^{\rm  RFIM}_N(\kappa,h,\gamma;J)\;,
\end{eqnarray*}
 is a well known consequence 
of  additivity and of the ergodic theorem, e.g. \cite{Vuill}.

Our main results can be summarized as follows:

\begin{theorem}
\label{th:0}
For all $h\in {\mathbb R}$, and all $\kappa, \beta >0$, the limit 
\begin{eqnarray}
\label{limite00}
p(\kappa,\beta,h)=\lim_{N\to\infty} \EE\,
p_N(\kappa,\beta,h;J)
%% \stackrel{N\to\infty}{\longrightarrow}
\end{eqnarray}
exists and 
\begin{eqnarray}
\label{limite01}
p_N(\kappa,\beta,h;J)\stackrel{N\to\infty}{\longrightarrow}
p(\kappa,\beta,h)
\end{eqnarray}
for a.e. $J$ and in the $L^p$-norm ($p \in [1,\infty)$). 
\end{theorem}
\begin{theorem}
\label{th:1}
For all $h\in {\mathbb R}$, there exist $\kappa_0(h)>0$ and $\beta_0(h)>0$ such that,
for $0\le \kappa\le \kappa_0(h)$ and $0\le \beta\le \beta_0(h)$
\begin{eqnarray}
\label{limite}
p(\kappa,\beta,h)
=
\inf_{0 \leq q \leq 1}
\left(p^{\rm  RFIM}(\kappa,h,\beta \sqrt{q})+\frac{\beta^2}4(1-q)^2\right)
\end{eqnarray} 
\end{theorem}

In the case $\kappa =0$, it is established in \cite{Gue00} that 
the infimum in (\ref{limite}) is achieved at a unique point. As explained 
below (see Section \ref{pf-unic}),
we take small enough $\kappa_0(h)$ 
 so that, 
for $0\le \kappa\le \kappa_0(h)$ and  $0\le \beta\le \beta_0(h)$,
the infimum in (\ref{limite}) is achieved at a unique $q$. In the following, we will always denote
by $\bar q=\bar q(\kappa,\beta,h)$ the value that realizes the infimum in (\ref{limite}).

\medskip

We emphasize that 
%%$\kappa_0(h), \beta_0(h)$ are taken small enough
$\kappa_0(h)$ is taken small enough
so that the RFIM is inside Dobrushin's uniqueness region for every realization
of the external fields $J_i, i\in {\mathbb Z}^d$. 
Let $\Med{\cdot}_{\infty,J}$ be the unique
 infinite volume Gibbs measure for the RFIM with $\gamma=\beta \sqrt{\bar q}$, 
depending on  $J_i, i\in {\mathbb Z}^d$.
If  $\preceq$ denotes the lexicographic order in ${\mathbb Z}^d$,
define
\begin{eqnarray*}
\Gamma=\Gamma(\kappa,\beta,h)=
\EE_{J_i,i\succeq 0} \left(\EE_{J'_0,J_i,i\prec 0}\log\Med{e^{\beta\sqrt{\bar q}
(J'_0-J_0)\sigma_0} }_{\infty,J}\right)^2,
\end{eqnarray*}
where $J'_0$ is an independent copy of $J_0$, $J'_0$ being independent of
$(J_i, i \in {\mathbb Z}^d)$. Here, the $\EE$-expectations 
are conditional, and the
subscripts of $\EE$
indicate on which variables the expectation is performed. Then, we have the following central limit theorem:
\begin{theorem}
\label{th:2}
For $0\le \kappa\le \kappa_0(h)$ and $0\le \beta\le \beta_0(h)$,
\begin{eqnarray}
\label{flutt}
\sqrt{|\vol|}\Big(p_N(\kappa,\beta,h;J)-\EE\,p_N(\kappa,\beta,h;J)\Big)
\stackrel{\rm law}{\longrightarrow}{\mathcal N}\left(0,\Gamma-\frac{\beta^2}2 \bar q^2\right).
\end{eqnarray}
\end{theorem}
One can check by expansion around $\kappa=\beta=0, h \neq 0$, that
the limiting variance $\Gamma-\frac{\beta^2}2 \bar q^2$ is strictly positive
for any fixed non-zero $h$ and small $\kappa, \beta$. This proves that 
the fluctuations of $p_N(\kappa,\beta,h;J)$ are truly of order
of the square root of the volume inverse 
in this region of the parameters. 
To match with the breakdown in the order of magnitude of fluctuations 
at zero external field,
we observe that  both $\bar q$ and $\Gamma$ vanish as $h \to 0$, and
also that they are equal to zero when $h=0$.

\section{Proofs}
\subsection{{\sl For small $\kappa$ the RFIM is inside the Dobrushin uniqueness region}
} \label{Dobr}

Let $i \neq k$ be lattice points. 
Under any (infinite volume) RFIM Gibbs measure  $\Med{\cdot}_{\infty,J}$
the law $\Med{\cdot | \eta}_{i,\infty,J}$
of $\sigma_i$ given $\sigma_j=\eta_j, j \neq i$, is
Bernoulli with parameter proportional to
$\exp\{ \sigma_i [{\cal H}_{i,k;\eta} + \kappa K(i-k)\eta_k]\}$
with
$${\cal H}_{i,k;\sigma}=
\kappa \sum_{j \neq i, k} K(i-j) \sigma_j + h + \gamma J_i.$$
Following e.g. Section 2 in \cite{Follmer}, we define the Dobrushin's
influence coefficient
$$
C_{ki}= \sup \left\{ \frac{1}{2} \|\Med{\cdot | \eta}_{i,\infty,J}-
\Med{\cdot | \eta'}_{i,\infty,J}\|_{\rm var} ; \eta=\eta' \;{\rm off \ }
k\right\}
$$
where $\|\cdot \|_{\rm var}$ is the variation norm. With a straightforward 
computation,
\begin{eqnarray*}
\| \Med{\cdot | \eta}_{i,\infty,J}-
\Med{\cdot | \eta'}_{i,\infty,J}\|_{\rm var}
 &=&
2| \Med{+| \eta}_{i,\infty,J}- \Med{+ | \eta'}_{i,\infty,J}|
\\
&\le&
\frac{ 2|\sinh 2 \kappa  K(i-k)|}{\cosh [2 {\cal H}_{i,k;\eta}]
+ \cosh [ 2 \kappa  K(i-k)]}
\end{eqnarray*}
 for such $\eta, \eta'$, and so
\begin{eqnarray}
\label{follmer}
C_{ki} \leq |\tanh  2 \kappa  K(i-k)| \leq  | 2 \kappa  K(i-k)|
\end{eqnarray}
Therefore, for $\kappa < \kappa_1=(2\sum_{i\neq 0} |K(i)|)^{-1}$
we derive from (\ref{follmer}) and (\ref{expdec}) that 
$a=\sup_i \sum_k 
\rho^{|i-k|}C_{ki}<1$
for some $\rho>1$,
which implies (see Section 2.3 in  \cite{Follmer}) that the Gibbs measure $ \Med{\cdot }_{\infty,J}$
is unique and has exponentially
decreasing correlations. More precisely,
for all local functions $f,g$, there is 
a finite constant   $C=C(a,f,g)$ such that for all $i\in \mathbb Z^d$
\begin{eqnarray}
  \label{eq:expdecr}
   |\Med{f;g \circ \theta_i}_{\infty,J}|\leq C \rho^{-|i|}\;,  
\end{eqnarray}
with $\theta_i$ the shift of vector $i$ and
$\Med{f;g}_{\infty,J}$ the covariance of $f,g$.

\subsection{{\sl Uniqueness of $\bar q$}} \label{pf-unic}

Introduce
$$
F(q,\kappa)=F_{\beta,h}(q,\kappa)= 
p^{\rm  RFIM}(\kappa,h,\beta \sqrt{q})+\frac{\beta^2}4(1-q)^2
$$
which is, in view of 
(\ref{eq:expdecr}), 
a smooth function of all its arguments if $\kappa < \kappa_1$
(e.g., Cor. 8.37 in \cite{Geo}).
\begin{proposition}\label{prop:1}
Assume %$\sum_{i\neq 0} K(i) \neq 0$ and 
$h \neq 0$. For all $\beta$, there exists
$\kappa_2=\kappa_2(\beta, h)>0$ such that 
$q \mapsto F(q,\kappa)$ has a unique minimizer on $[0,1]$
for $\kappa < \kappa_2$.
\end{proposition}

By \cite{Gue00}, we know that the minimizer $q_0$ of 
$F(q,0)$ is unique, strictly positive and that\footnote{ \cite{Gue00}, p.166}
%%\footnote{ \cite{Gue00}, p.66}
\begin{eqnarray}
\label{implicit}
\frac{\partial^2}{\partial^2 q}F(q_0,0)>0.
\end{eqnarray}
 By (\ref{eq:expdecr}),
the function  $F$ is
continuous in $\kappa$ uniformly in $q \in  (0,1]$. Hence, 
\begin{eqnarray} \label{eq:cvunif}
\forall \delta>0 \;\; \exists \tilde \kappa:\; \forall \kappa < \tilde\kappa\quad
{\rm argmin}_{[0,1]}\; F(\cdot, \kappa) \subset (q_0-\delta, q_0+\delta)
\end{eqnarray}
Again by (\ref{eq:expdecr}),
the  function  $F$ is
$C^2$ in a neighborhood of $(q_0,0)$. 
By the implicit function theorem for the equation
$$
\frac{\partial}{\partial q} F(q,\kappa)=0\;
$$
(which applies thanks to condition (\ref{implicit}))
we then derive that there exist neighborhoods $U$ of 
$q_0$ and $V$ of $\kappa=0$ and a function $\bar q: V \mapsto U$
such that, for $q \in V, \kappa \in U$,
the above equation is equivalent to $q =\bar q (\kappa)$. 
With $\delta$ small enough so that $(q_0-\delta, q_0+\delta)
\subset V$, we choose now $\kappa_2<\tilde \kappa$ (with $\tilde \kappa$
from (\ref{eq:cvunif})) such that
$(-\kappa_2,\kappa_2) \subset U$. 
Then, for $\kappa< \kappa_2$, $F(\cdot, \kappa)$ has a unique
minimum at $q=\bar q(\kappa)$.

As for the case $h=0$, one can prove similarly 
the following results, that we mention for comparison but will neither prove
nor use. 
\begin{proposition}
For all $\beta\ne1$, there exists
$\kappa_3=\kappa_3(\beta)>0$ such that 
$q \mapsto F_{\beta,0}(q,\kappa)$ has a unique minimizer on $[0,1]$
for $\kappa < \kappa_3$. If $\beta<1$, the minimizer is $q=0$.
\end{proposition}

{\bf Remark} The restriction $\beta\ne1$ is due to the fact that
for $\beta=\beta_c=1$ (the critical point of the Sherrington-Kirkpatrick model)
one has $q_0(\beta_c,0)=0$ and, in contrast with (\ref{implicit}),
$$
\frac{\partial^2}{\partial^2 q}F_{\beta_c,0}(0,0)=0.
$$

\subsection{{\sl Proof of Theorem \ref{th:0}} }

{\em Proof of (\ref{limite00})} The proof is standard, we just sketch 
the main steps. 
Consider the  box $\Lambda_{mN}$, with $m,N\in\mathbb N$, and partition it into sub-boxes $\vol^{(\ell)}$, 
$\ell=1,\cdots,m^d $,
congruent to $\vol$. Moreover, let $\tilde Z_{mN}$ be the partition function of the system where, with respect to
(\ref{H}), $H_{mN}^{SK}(\sigma;J)$ is replaced by
$$
-\sum_{\ell=1}^{m^d}\frac1{\sqrt{2|\vol|}}\sum_{i,j\in \vol^{(\ell)}}J_{ij}^{(\ell)}\sigma_i\sigma_j,
$$
the $J^{(\ell)}_{ij}$ being $m^d$ independent families of standard Gaussian variables. 
Note that, in this system, the different sub-boxes interact only through the Ising potential $K(.)$.
Then, following the ideas of \cite{GueTon02pres}, it is easy to prove that
\begin{eqnarray*}
p_{mN}(\kappa,\beta,h)\ge \frac{1}{|\Lambda_{mN}|}\EE\,\log \tilde Z_{mN}(\kappa,\beta,h;J^{(\ell)}).
\end{eqnarray*}
Since the Ising potential is summable, the interaction among the different sub-boxes due to the potential $K$ 
grows at most proportionally to $\kappa$ and 
to the the total surface, $d \,m^{d} N^{d-1}$. As a consequence, one has
the approximate monotonicity
\begin{eqnarray*}
p_{mN}(\kappa,\beta,h)\ge p_{N}(\kappa,\beta,h)-\kappa\frac{C}{N},
\end{eqnarray*}
for some constant $C$ depending on the potential $K(.)$.
From this, it is a standard fact to deduce that the sequence $\left\{p_{mN}\right\}_N$
has a limit when $N\to\infty$, that it does not depend on $m$, and that it coincides with
$\lim_N p_N$.

% (To  be written: Follow the subadditive 
% proof of Guerra-Toninelli 
% for SK, to compare pressure for Ising on box $\Lambda_n$ with 
% SK there,  and Ising on box $\Lambda_n$ with two uncoupled
% SK on a partition of  $\Lambda_n$. Apply this to a 
%  partition of  $\Lambda_{mn}$ into $m^d$ boxes like  $\Lambda_n$.
% At last, neglect the Ising interaction between the latter small boxes,
% to get approximate subadditivity.)

{\sl Proof of (\ref{limite01})} 
The almost sure convergence  is standard and follows from 
exponential self-averaging of the pressure (see for instance Proposition 2.18 of \cite{ledoux}),
which in the present case reads
$$
\mathbb P\,\left(\left|p_N(\kappa,\beta,h;J)-\EE\,p_N(\kappa,\beta,h;J)\right|\ge u\right)
\le D_1\; e^{-D_2(\beta)|\vol|u^2},
$$
together with Borel-Cantelli's lemma.
The $L^p$-convergence comes from uniform integrability, which
again follows from exponential concentration.

\subsection{{\sl Proof of Theorem \ref{th:1}}}

%The proof is based on a method (``quadratic replica coupling'') introduced in 
%\cite{GueTon02quad}.
For $0\le t\le1$ and any $q\ge0$, define the interpolating partition function
\begin{eqnarray}
\label{interpola}
Z(t)= 
\sum_{\sigma\in\{-1,+1\}^{\Lambda_N}}\exp\Big({- H^{(t)}(\sigma)}\Big)
\end{eqnarray}
$$
H^{(t)}(\sigma)= \kappa H^I_N(\sigma)
+\beta \sqrt t H^{SK}_N(\sigma;J)
-\sum_{i\in \Lambda_N}\sigma_i
\left(h+\beta \sqrt{q(1-t)}J_i\right)
$$
with the properties 
\begin{eqnarray*}
&& Z(0)=Z_N^{\rm  RFIM}(\kappa,h,\beta \sqrt{q};J)\\
&& Z(1)=Z_N(\kappa,\beta,h;J).
\end{eqnarray*}
The $t$-derivative of the corresponding pressure 
$$
p_N(t)=\frac1{|\vol|}\EE\,\log Z(t)
$$
is easily computed: We denote by  $ \Med{\cdot}_t$ the Gibbs measure
associated to $H^{(t)}$, by $ \Med{\cdot}_t^{\otimes 2}$ its tensor product
acting on a pair $(\sigma^1,\sigma^2)\in\{-1,+1\}^{\Lambda_N}\times
\{-1,+1\}^{\Lambda_N}$, by $q_{12}=|\Lambda_N|^{-1}
\sum_{i\in \Lambda_N}\sigma^1_i \sigma^2_i$ the overlap between configurations $\sigma^1$ and
$\sigma^2$, and we 
get by the Gaussian integration by parts formula $\EE\,J_{i}F(J_i)=
\EE\,F'(J_i)$,
\begin{eqnarray}
\frac d{dt}p_N(t)
&=& \nonumber 
\frac{\beta}{2 |\Lambda_N|} \EE
\left\{ \frac{1}{\sqrt{ 2 |\Lambda_N|t}} 
\sum_{i,j \in \Lambda_N} J_{ij}
\Med{\sigma_i \sigma_j}_t
-\frac{
\sqrt q}{\sqrt{1-t}}
\sum_{i\in \Lambda_N}
  J_i\Med{\sigma_i}_t
\right\}\\
&
\stackrel{\rm int.by\,parts}{=}
&  \nonumber 
\frac{\beta}{2 |\Lambda_N|} \EE
\left\{ \frac{\beta}{{ 2 |\Lambda_N|}} 
\sum_{i,j \in \Lambda_N} (1-
\Med{\sigma^1_i \sigma^1_j\sigma^2_i \sigma^2_j}_t^{\otimes 2})
-\beta q
\sum_{i\in \Lambda_N} (1-
 \Med{\sigma^1_i \sigma^2_i}_t^{\otimes 2})
\right\}\\
&=&  \label{i2}
\frac{\beta^2}4(1-q)^2-\frac{\beta^2}4 
\EE\,\Med{(q_{12}-q)^2}_t^{\otimes 2}
\end{eqnarray}
so that, integrating in $t$ between $0$ and $1$, taking the $N\to\infty$ limit and optimizing on 
$q$, we have the ``first half'' of Eq. (\ref{limite}):
\begin{eqnarray}
\label{half}
p(\kappa,\beta,h)
\le
\inf_{q\ge0}\left(p^{\rm  RFIM}(\kappa,h,\beta \sqrt{q})+\frac{\beta^2}4(1-q)^2\right).
\end{eqnarray}
Note that the inequality holds for any values of $\beta$ and $\kappa$. 

\medskip

Next, consider a system of two coupled replicas: For  $\lambda>0$ let 
\begin{eqnarray*}
Z^{(2)}(t,\lambda)=\sum_{\sigma^1,\sigma^2\in\{-1,+1\}^{\Lambda_N}}\exp
\left(-H^{(t)}(\sigma^1)-H^{(t)}(\sigma^2)+\frac{\beta^2}2|\vol|\lambda(q_{12}-q)^2
\right)
\end{eqnarray*}
and denote by  
$ \Med{\cdot}_{t,\lambda}$ the Gibbs measure on $\{-1,+1\}^{\Lambda_N}\times
\{-1,+1\}^{\Lambda_N}$
associated to this partition function. (Of course, 
$ \Med{\cdot}_{t,0}=   \Med{\cdot}_t^{\otimes 2}$ and $Z^{(2)}(t,0)=
Z(t)^2$.) With 
$$p^{(2)}_N(t,\lambda)=\frac1{2|\vol|}\EE\,\log Z^{(2)}(t,\lambda)$$
we have by a  computation similar to (\ref{i2}) and using symmetry
in $\sigma^1, \sigma^2$,
\begin{eqnarray*}
&\qquad \qquad\frac{\displaystyle d}{\displaystyle dt}&
p^{(2)}_N
(t,\lambda_0-t) = \qquad\qquad\\
&=&
\frac{\beta}{4 |\Lambda_N|} \EE
\left\{ \frac{2}{\sqrt{ 2 |\Lambda_N|t}} 
\sum_{i,j \in \Lambda_N} J_{ij}
\Med{\sigma^1_i \sigma^1_j}_{t,\lambda_0-t}
-\frac{2
\sqrt q}{\sqrt{1-t}}
\sum_{i\in \Lambda_N}
  J_i\Med{\sigma^1_i}_{t,\lambda_0-t} \right. \\
&& \left. \qquad \qquad
- \beta |\Lambda_N| \Med{(q_{12}-q)^2
}_{t,\lambda_0-t}
\right\}\\
&
\stackrel{\rm int.by\,parts}{=}
&  \nonumber 
\frac{\beta}{4 |\Lambda_N|} \EE
\left\{ \frac{\beta}{{  |\Lambda_N|}} 
\sum_{i,j \in \Lambda_N}  (1+
\Med{\sigma^1_i \sigma^1_j\sigma^2_i \sigma^2_j}_{t,\lambda_0-t}
-2
\Med{\sigma^1_i \sigma^1_j\sigma^3_i \sigma^3_j}_{t,\lambda_0-t}^{\otimes 2})
\right. \\
&& \left. \qquad 
-2\beta q
\sum_{i\in \Lambda_N} (1+ \Med{\sigma^1_i \sigma^2_i}_{t,\lambda_0-t}
-2
 \Med{\sigma^1_i \sigma^3_i}_{t,\lambda_0-t}^{\otimes 2})
- \beta |\Lambda_N| \Med{(q_{12}-q)^2
}_{t,\lambda_0-t}
\right\}\\
&=&
\frac{\beta^2}4(1-q)^2-\frac{\beta^2}2
\EE\,\Med{(q_{13}-q)^2}_{t,\lambda_0-t}^{\otimes 2}
\end{eqnarray*}
so that, integrating,
\begin{eqnarray} \label{01}
%%\frac1{2|\vol|}\EE\,\log Z^{(2)}(t,\lambda)&\le& 
p^{(2)}_N(t,\lambda) \leq
\frac{\beta^2}4(1-q)^2t+\frac1{2|\vol|}
\EE\,\log Z^{(2)}(0,t+\lambda).
%&\equiv& p_N^{RFIM}(t)+\frac1{2|\vol|}\EE\,\log \left(\frac{Z^{(2)}(0,t+\lambda)}
%{Z^{(2)}(0,0)}\right)
\end{eqnarray}
Now, setting $u_N(t)= p_N^{\rm  RFIM}(\kappa,h,\beta \sqrt{q})+
\frac{\beta^2}4(1-q)^2t-p_N(t)$, which is a non negative function by  
(\ref{i2}), and using convexity of the pressure with respect to $\lambda$
and the identity $p_N^{(2)}(t,0)=p_N(t)$,
we obtain for any $\lambda>0$
\begin{eqnarray}
\label{gronw}
%% 0&{\leq}& \frac{\beta^2}4\EE\,\Med{(q_{12}-q)^2}_t
%% \stackrel{(\ref{i2})}{=} 
%%\frac d{dt}\left(p_N^{\rm  RFIM}(\kappa,h,\beta \sqrt{q})
%%+\frac{\beta^2}4(1-q)^2t-p_N(t)\right)\\\nonumber 
\frac d{dt}u_N(t)
&\stackrel{(\ref{i2})}{=}& 
 \frac {\partial }{\partial \lambda}p_N^{(2)}(t, 0) \nonumber \\
&\stackrel{\rm convexity}{\leq} &
\frac{p^{(2)}_N(t,\lambda)-p_N(t)}\lambda\\\nonumber
&\stackrel{(\ref{01})}{\leq} &
\frac1\lambda\left[
%%{p_N^{\rm  RFIM}(\kappa,h,\beta \sqrt{q})+\frac{\beta^2}4(1-q)^2t-p_N(t)
u_N(t)
+\frac1{2|\vol|}\EE\,\log \frac{Z^{(2)}(0,t+\lambda)}
{Z^{(2)}(0,0)}\right].
\end{eqnarray}
%By  Gronwall's lemma, and 
Since $Z^{(2)}$ is increasing in $\lambda$, Eq. (\ref{gronw}) implies
\begin{eqnarray*}
\frac d{dt}\log\left[u_N(t)+\frac1{2|\vol|}\EE\,\log \frac{Z^{(2)}(0,1+\lambda)}{Z^{(2)}(0,0)}\right]\le\frac1\lambda,
\end{eqnarray*}
which, recalling that $u_N(0)=0$, can be immediately integrated to give
\begin{eqnarray}
\label{gronw2}
u_N(t) \leq \left(e^{t/\lambda}-1\right)\times\frac1{2|\vol|}\EE\,\log \frac{Z^{(2)}(0,1+\lambda)}
{Z^{(2)}(0,0)}
%\lambda(e^{\lambda t}-1) \times
%\frac1{2|\vol|}\EE\,\log \frac{Z^{(2)}(0,t+\lambda)}
%{Z^{(2)}(0,0)}
\end{eqnarray}
and Eq. (\ref{limite}) follows  if we can prove that
\begin{eqnarray}
\label{eq:provanda}
\lim_{N\to\infty}\frac1{2|\vol|}\EE\,\log \frac{Z^{(2)}(0,\lambda_0)}
{Z^{(2)}(0,0)}=0
\end{eqnarray}
for some $\lambda_0>1$, if $q$ is chosen properly. 
%%($\lambda_0=\lambda+t>1$ 
(Note that $\lambda_0=\lambda+1>1$ 
is required so that one can take
$t$ up to $1$ and still have $\lambda>0$, which is needed in (\ref{gronw}).)
%%%%).

Define for any $\mu\in \mathbb R$, $q\ge0$
\begin{eqnarray*}
\alpha_N(\mu;J)=\frac1{2|\Lambda_N|} \log \Med{e^{\mu|\Lambda_N|(q_{12}-q)}}^{\otimes 2}_{\kappa,\beta \sqrt{q},N}
\end{eqnarray*}
and $\alpha_N(\mu)=\EE\,\alpha_N(\mu;J)$. We denote by
$\Med{.}^{\otimes 2}_{\kappa,\beta \sqrt{q},N}$ 
 the Gibbs measure for two
replicas of the RFIM with parameters $\kappa$ and 
$\gamma=\beta \sqrt{q}$ and volume 
$\Lambda_N$. Later we will also use the notation 
$\Med{.}_{\kappa,\beta \sqrt{q},N}^{(\mu)}$
for 
\begin{eqnarray*}
\Med{A}_{\kappa,\beta \sqrt{q},N}^{(\mu)}=
\frac{\Med{A\, e^{\mu|\Lambda_N|(q_{12}-q)}}^{\otimes 2}_{\kappa,\beta \sqrt{q},N}}
{\Med{e^{\mu|\Lambda_N|(q_{12}-q)}}^{\otimes 2}_{\kappa,\beta \sqrt{q},N}}.
\end{eqnarray*}

Let  $\bar q_N=\bar q_N(\beta,\kappa,h)$ be the value which 
 minimizes
$$
p_N^{\rm  RFIM}(\kappa,h,\beta \sqrt{q})+\frac{\beta^2}4(1-q)^2
$$
with respect to $q$, cf (\ref{limite}). Clearly, $\bar q_N$ satisfies 
the ``self-consistent equation''
\begin{eqnarray}
\label{inf}
q=\EE\Med{q_{12}}^{\otimes 2}_{\kappa,\beta \sqrt{ q},N}
= \frac{\sum_{i\in\Lambda_N}\EE\Med{\sigma_i}^2_{\kappa,\beta \sqrt{ q},N}}{|\Lambda_N|}.
%%\lim_{N\to\infty}\EE\,\Med{\sigma_0}^2_{\kappa,\beta \sqrt{q},N} 
\end{eqnarray}
An analysis analogous to the one of Section \ref{pf-unic} shows that 
the solution of (\ref{inf}) is unique for $\kappa$ small enough, for $N$ sufficiently large,
and that $\bar q_N\to\bar q$ for $N\to\infty$.
%%%\footnote{ (check *******. I think it is ok).}
%%% yes !
A Taylor expansion around $\mu=0$ gives immediately
\begin{eqnarray*}
\alpha_N(\mu;J)=\alpha_N(0;J)+\mu\alpha'_N(0;J)+\int_0^\mu dy \int_0^y du\, \alpha''_N(u;J),
\end{eqnarray*}
where
\begin{eqnarray}
\alpha_N(0;J)&=&0\\
\label{deriv1}
\alpha'_N(0;J)&=&\frac12\left(\Med{q_{12}}^{\otimes 2}_{\kappa,\beta \sqrt{\bar q_N},N}-
\EE\,\Med{q_{12}}^{\otimes 2}_{\kappa,\beta \sqrt{\bar q_N},N}\right)\\
\alpha''_N(u;J)&=&\frac{|\Lambda_N|}2 \left(\Med{q_{12}^2}^{(u)}_{\kappa,\beta \sqrt{\bar q_N},N}-
(\Med{q_{12}}^{(u)}_{\kappa,\beta \sqrt{\bar q_N},N})^2\right)\\
&=&\frac1{2|\Lambda_N|}\sum_{i,j\in\Lambda_N}\Med{(\sigma^1_i\sigma^2_i-\Med{\sigma^1_i\sigma^2_i})
(\sigma^1_j\sigma^2_j-\Med{\sigma^1_j\sigma^2_j})}^{(u)}_{\kappa,\beta \sqrt{\bar q_N},N}.
\end{eqnarray}
We can view $\Med{.}_{\kappa,\beta \sqrt{\bar q_N},N}^{(u)}$ as the Gibbs 
measure of a system with $2|\Lambda_N|$ spins, 
with exponentially decaying pair interactions. Since  $\kappa$ is small,
taking $\mu$ itself sufficiently small, we keep this system 
inside the Dobrushin uniqueness region.
Using exponential decay of correlations as in Section \ref{Dobr}, we 
obtain that $\alpha''_N(u;J)$ is bounded above by a constant, uniformly in $N$,$J$ and in $u \in [0, \mu]$,
so that
\begin{eqnarray}
\label{tailor}
\alpha_N(\mu;J)\le \mu\alpha'_N(0;J)+C \mu^2
\end{eqnarray}
for $\kappa<\kappa_0(h)$. Note that this bound
holds for any $\mu$, since it does for small $\mu$ and the function
$\alpha$ can grow at most linearly at infinity.

With this in hand, 
we go back to proving (\ref{eq:provanda}).
After a Gaussian transformation
\begin{eqnarray*}
\frac{Z^{(2)}(0,\lambda_0)}
{Z^{(2)}(0,0)}&=&
\int dz\sqrt{\frac{|\Lambda_N|}{2\pi}} 
e^{-|\Lambda_N| \frac{z^2}2}\Med{e^{\beta z \sqrt{\lambda_0}|\Lambda_N|(q_{12}-\bar q_N)}}^{\otimes 2}
_{\kappa,\beta \sqrt{\bar q_N},N}\\
&=&  \int dz\sqrt{\frac{|\Lambda_N|}{2\pi}} 
e^{|\Lambda_N|(- z^2/2+2\alpha_N(\beta z \sqrt{\lambda_0};J))}.
\end{eqnarray*}
Equations (\ref{tailor}) and (\ref{deriv1}) imply that, for $\beta<\beta_0(h)=1/\sqrt{4C\lambda_0}$
\begin{eqnarray}
\nonumber
\frac1{2|\vol|}\EE\,\log \frac{Z^{(2)}(0,\lambda_0)}
{Z^{(2)}(0,0)}&\le&\frac1{2|\vol|}\EE\,\log \int dz \sqrt{\frac{|\Lambda_N|}{2\pi}}
 e^{|\Lambda_N|\left(-z^2/2(1-4C\beta^2\lambda_0)+2z\beta \sqrt{\lambda_0}\alpha'_N(0;J)\right)}\\
\nonumber
&\le& C'\EE\,(\alpha'_N(0;J))^2\\
&=&
\label{limitanda}
\frac{C'}{4|\Lambda_N|^2}\sum_{i,j\in\Lambda_N}\EE((\Med{\sigma_i}^2-\EE\Med{\sigma_i}^2)
(\Med{\sigma_j}^2-\EE\Med{\sigma_j}^2)),
\end{eqnarray}
where for simplicity we have written $\med.$ for $\med.^{\otimes 2}_{\kappa,\beta \sqrt{\bar q_N},N}$.
We now show that 
the last expression is of order $1/|\vol|$. 
%%Indeed, take two distinct sites $i,j$ and
%%consider the $d$-dimensional sphere $S_{ij}$ of radius $|i-j|/2$ centered at $i$, so that both $i$ and $j$ 
%%have distance at least $|i-j|/2$ from any point of the surface of $S_{ij}$. If the interactions
%%between the interior and the exterior of $S_{ij}$ were absent, the corresponding term 
%in (\ref{limitanda}) would be zero. Then, using the exponential decay of correlation (\ref{eq:expdecr}) it
%% is 
%%%%%%%%%%%%%%%%%%%%%%%%%%%%%%%%%%%%%%%%%%%%immediate 
%% to see that
Indeed, take two distinct sites $i,j$ and
consider the $d$-dimensional ball $B_{ij}$ of radius  $|i-j|/2$ 
centered at $i$. If $\Med{\sigma_i}$ were depending only on
$(J_k, k \in B_{ij})$ and $\Med{\sigma_j}$ on $(J_k, k \in B_{ij}^c)$ only,
the corresponding term in (\ref{limitanda}) would be zero by independence.
But in the Dobrushin region, we can approximate $\Med{\sigma_i}$
by the expectation of $\sigma_i$ for the finite volume RFIM
on $B_{ij}$ with an error which is exponentially small in the radius
 $|i-j|/2$, uniformly in the $J_k$'s. Doing similarly with 
 $\Med{\sigma_j}$, we conclude that
$$
\EE((\Med{\sigma_i}^2-\EE\Med{\sigma_i}^2)
(\Med{\sigma_j}^2-\EE\Med{\sigma_j}^2))\le C' \rho^{-C''|i-j|}
$$
for suitable constants $C',C''>0$. 
This, together with (\ref{limitanda}), immediately implies that
\begin{eqnarray*}
\frac1{2|\vol|}\EE\,\log \frac{Z^{(2)}(0,\lambda_0)}{Z^{(2)}(0,0)}\le \frac{c}{|\vol|}.
\end{eqnarray*}
At this point, recalling Eq. (\ref{gronw2}), one finds 
\begin{eqnarray*}
p_N(\kappa,\beta,h)\ge p_N^{RFIM}(\kappa,h,\beta\sqrt{\bar q_N})+\frac{\beta^2}4(1-\bar q_N)^2-\left(
e^{1/(\lambda_0-1)}-1\right)\frac c{|\vol|},
\end{eqnarray*}
which, together with (\ref{half}), proves the 
convergence in average of the pressure to the expression  (\ref{limite}).
Moreover, from Eqs. (\ref{gronw}) and (\ref{i2}) one deduces that
\begin{eqnarray}
\label{1N}
\limsup_{N\to\infty} \sup_{0\le t \le1}
{|\vol|}\EE\,\Med{(q_{12}-\bar q_N)^2}^{\otimes 2}_t <\infty\;,
\end{eqnarray}
which will be needed in next section, where we deal with pressure fluctuations.

\subsection{{\sl Proof of Theorem \ref{th:2}}}

We follow the strategy which was introduced in \cite{GueTon02tcl}, adapted to the present case where
short-range interactions are also present.
Let 
\begin{eqnarray*}
\hat f_N(t)=\sqrt{|\Lambda_N|}\left(\frac{\log Z(t)}{|\Lambda_N|}-p_N(t) 
\right),
\end{eqnarray*}
where $Z(t), p_N(t)$ were defined in Eqs. (\ref{interpola})-(\ref{i2}) and it is understood that 
$q$ is taken to be $\bar q_N=\bar q_N(\beta,\kappa,h)$ as in Eq. (\ref{inf}), to be distinguished 
from $\bar q=\lim \bar q_N$.

We will prove that
\begin{eqnarray}
\label{intermedia}
\lim_{N\to\infty} \EE\,e^{iu \hat f_N(t)}=\exp\left(-\frac{u^2}2(\Gamma-\frac{\beta^2}2t \bar q^2)
\right)
\end{eqnarray}
for any $u\in \mathbb R$, $0\le t\le 1$, from which Eq. (\ref{flutt}) follows for $t=1$.

By means of integrations by parts one finds
\begin{eqnarray*}
&&\partial_t\EE\,e^{i u \hat f_N(t)}=iu\EE\,e^{i u \hat f_N(t)}\frac d{dt} \hat f_N(t)\\
&&=\frac{\beta^2}4u^2\bar q_N^2\EE\,e^{i u \hat f_N(t)}-
\frac{\beta^2}4u^2\EE\,e^{i u \hat f_N(t)} \Med{(q_{12}-\bar q_N)^2}^{\otimes 2}_t\\
&&-i \frac{\beta^2}4 u \sqrt{|\vol|}\EE\,e^{i u \hat f_N(t)}\left(\Med{(q_{12}-\bar q_N)^2}^{\otimes 2}_t-
\EE\Med{(q_{12}-\bar q_N)^2}^{\otimes 2}_t\right)
\end{eqnarray*}
and, using Eq. (\ref{1N}), 
\begin{eqnarray*}
\partial_t\EE\,e^{i u \hat f_N(t)}=\frac{\beta^2u^2\bar q^2}4\EE\,e^{i u \hat f_N(t)}+o(1).
\end{eqnarray*}
Integrating in $t$, Eq. (\ref{intermedia}) is then proven provided that we show that
\begin{eqnarray}
\label{provided}
\lim_{N\to\infty}\EE\,e^{iu \hat f_N(0)}=\exp\left(-\frac{u^2}2\Gamma
\right),
\end{eqnarray}
i.e., a central limit theorem for pressure fluctuations of the RFIM at high temperature.

To this purpose, we employ the central limit theorem for martingales \cite{clt_marti}, 
which we recall for convenience.
Let $(\Omega,{\cal A},P)$ be a probability space,   
${\cal F}^{(n)}=\left\{{\cal F}^{n}_k\right\}_{0\le k\le n}$ 
a filtration of $\cal A$, for $n\in \mathbb N$, such  that ${\cal F}^{n}_0=\left\{\emptyset,\Omega\right\}$,
and $\xi^{(n)}=\left\{\xi_{n,k}\right\}_{1\le k\le n}$ a sequence of random variables adapted to ${\cal F}^{(n)}$.
We denote by $E^{n,k}$ (respectively $P^{n,k}$) the expectation (respectively the probability) 
conditioned to ${\cal F}^{n}_k$ and by $V^{n,k}$ the conditional variance
\begin{eqnarray*}
V^{n,k}(X)=E^{n,k}(X^2)-\left(E^{n,k} X\right)^2,
\end{eqnarray*}
of a random variable $X$.
We say that the triangular array $\left\{\xi_{n,k}\right\}_{n>0,1\le k\le n}$ is asymptotically negligible if
for any $\varepsilon>0$
\begin{eqnarray}
\label{asneg}
\sum_{k=1}^n P^{n,k-1}\left(|\xi_{n,k}|\ge\varepsilon\right)\stackrel{\rm P}{\longrightarrow}0
\end{eqnarray}
when $n\to\infty$.
Then, the following holds:
\begin{theorem}
\label{marty}
Let $\left\{\xi_{n,k}\right\}_{n>0,1\le k\le n}$ be an asymptotically negligible triangular array
of square integrable random variables, and assume that for some $\Gamma>0$, 
\begin{eqnarray}
\label{condiz1}
\sum_{k=1}^n E^{n,k-1} (\xi_{n,k})\stackrel{\rm P}{\longrightarrow}0
\end{eqnarray}
and
\begin{eqnarray}
\label{condiz2}
\sum_{k=1}^n V^{n,k-1}(\xi_{n,k})\stackrel{\rm P}{\longrightarrow}\Gamma,
\end{eqnarray}
for $n\to\infty$.
Then, 
\begin{eqnarray*}
\sum_{k=1}^n\xi_{n,k}\stackrel{\rm law}{\longrightarrow}\mathcal N (0,\Gamma).
\end{eqnarray*}
\end{theorem}

To simplify notations in our case, let $|\vol|=n$, $h+\beta \sqrt {\bar q_N} J_i=h_i$  
and 
$$p_{n,\uh}=p^{RFIM}_N(\kappa,h,\beta \sqrt{\bar q_N};J).$$

%%%%%%%%%%%% corrections d6

Introducing the usual lexicographic ordering of the 
sites in $\vol$, we define
$
{\cal F}^{n}_k
$, for $k\in \vol$,
as the $\sigma$-algebra generated by the random fields $J_i$ for $i\in \vol,i\preceq k$, and
\begin{eqnarray*}
\xi_{n,k}=\sqrt n(E^{n,k}p_{n,\uh}-E^{n,k-1}p_{n,\uh})
\end{eqnarray*}
so that
\begin{eqnarray}
\label{riscalata}
{\sqrt{|\vol|}}({p^{RFIM}_N(\kappa,h,\beta \sqrt{\bar q_N};J)-\EE\,p^{RFIM}_N(\kappa,h,\beta \sqrt{\bar q_N};J)})
=\sum_{k\in\vol}\xi_{n,k}.
\end{eqnarray}
Of course, $E^{n,k}(.)=\EE_{J_{\ell},\ell\in \vol,\ell\succ k}(.)$, and by convention
${\cal F}^{n}_{k}=\{\emptyset,\Omega\}$
if $k$ precedes the first site in $\vol$.
One can rewrite
\begin{eqnarray*}
\xi_{n,k}=-\frac1{\sqrt n}\EE_{J'_k,J_\ell,\ell\succ k}\log \Med{e^{({h'}_k-h_k)\sigma_k}}_{n,\uh}
\end{eqnarray*}
where $J'_k$ is an independent copy of $J_k$
-- independent of $(J_l, l \in {\mathbb Z}^d)$ --,
 and $h'_k=h+\beta \sqrt{\bar q_N}J'_k$, so that
\begin{eqnarray*}
\left|\xi_{n,k}\right|\le n^{-1/2}\EE_{J'_k}|h'_k-h_k|\le Cn^{-1/2}(|h_k|+C)
\end{eqnarray*}
for some constant $C$ and
\begin{eqnarray*}
E^{n,k-1} \left|\xi_{n,k}\right|^3\le C' n^{-3/2}.
\end{eqnarray*}
This implies asymptotic negligibility (\ref{asneg}), since
\begin{eqnarray*}
\sum_{k\in \vol}P^{n,k-1}\left(|\xi_{n,k}|\ge\varepsilon\right)\le \frac1{\varepsilon^3}
\sum_{k\in\vol} E^{n,k-1}\left(|\xi_{n,k}|^3\right)\le \frac{C'}{\varepsilon^3\sqrt n}.
\end{eqnarray*}
In order to apply Theorem \ref{marty}, we have to check conditions (\ref{condiz1}) and (\ref{condiz2}). The 
first one is evident, since
\begin{eqnarray*}
E^{n,k-1}\xi_{n,k}=0
\end{eqnarray*}
identically.
As for the second, notice that
\begin{eqnarray}
\label{var}
V^{n,k-1}(\xi_{n,k})=\frac1n\EE_{J_k}\left(\EE_{J'_k,J_\ell,\ell\succ k}\log \Med{
e^{(h'_k-h_k)\sigma_k}}_{n,\uh}
\right)^2.
\end{eqnarray}
Let $k$ correspond to a site ``in the bulk'' of $\vol$, 
i.e., assume that the distance between $k$ and the boundary of $\vol$ is larger than,
say, $n^{1/(2d)}$. In this case,   we will write $k\in B_N$.
% \begin{eqnarray}
% \lim_{n\to\infty}\sup_{\uh, k}||\Med{\{\sigma_k\in.\}}_{n,\uh}-\Med{\{\sigma_k\in.\}}_{\infty,\uh}||\equiv 
% \lim_{n\to\infty}\epsilon_n=0
% \end{eqnarray}
We want to replace  $\Med{.}_{n,\uh}$  in (\ref{var}) with the unique infinite-volume Gibbs measure 
$\Med{.}_{\infty,{\uh}}$.
To this purpose, note preliminarily that
\begin{eqnarray*}
\left|\log\frac{\Med{e^{( h'_k-h_k)\sigma_k}}_{n,\uh}}
{\Med{e^{( h'_k- h_k)\sigma_k}}_{\infty,{\uh}}}\right|\le 2| h'_k- h_k|.
\end{eqnarray*} 
Moreover, thanks to Dobrushin's theorem,
\begin{eqnarray} \label{palezo}
\frac{\Med{e^{(h'_k-h_k)\sigma_k}}_{n,\uh}}
{\Med{e^{(h'_k-h_k)\sigma_k}}_{\infty,\uh}}=1+\frac{e^{(h'_k-h_k)}-e^{-(h'_k-h_k)}}{
\Med{e^{( h'_k- h_k)\sigma_k}}_{\infty,\uh}}\delta_{n,\uh}
\end{eqnarray}
with some $\delta_{n,\uh}$ such that
$$
\lim_{n\to\infty}\sup_{k\in B_N}\sup_{\uh}|\delta_{n,\uh}|\equiv \lim_{n\to\infty}\epsilon_n=0.
$$ 
Denoting by $A_{n,k}$ the event 
$$
A_{n,k}=\left\{|h'_k-h_k|\le\frac12\log\left(\frac1{2\epsilon_n}\right)
\right\},
$$
and using the fact that $|\log(1+x)|\le D|x|$ for $|x|\le1/2$ for some finite constant $D$,
one can write
\begin{eqnarray*}
\left|\log\frac{\Med{e^{(h'_k-h_k)\sigma_k}}_{n,\uh}}
{\Med{e^{(h'_k-h_k)\sigma_k}}_{\infty,\uh}}\right|
\le D\epsilon_n e^{2|h'_k-h_k|}+2|h'_k-h_k|
1_{A_{n,k}^C}.
\end{eqnarray*}
Therefore one has
\begin{eqnarray*}
V^{n,k-1}(\xi_{n,k})&=&\frac1n\EE_{J_k}\left(\EE_{J'_k,J_\ell,\ell\succ k}\log \Med{
e^{(h'_k-h_k)\sigma_k}}_{\infty,\uh}
\right)^2+\frac1no(1)\\
&=&\frac1n\left(\phi(\theta_{-k} \uh)+o(1)\right),
\end{eqnarray*}
where
%$\EE\,1_{A_{n,k}^C}=o(1)$ 
$o(1)\to0$ for $n\to\infty$ uniformly in $k\in B_N$,
 $\theta_k$ is the shift of vector $k$ and
\begin{eqnarray*}
\phi({\uh})=\EE_{J_0}\left(\EE_{J'_0,J_\ell,\ell\succ 0}\log\Med{
e^{(h'_0-h_0)\sigma_0}}_{\infty,\uh}
\right)^2,
\end{eqnarray*}
the subscript $_0$ referring of course to the origin of the lattice $\mathbb Z^d$.
Note that there is a residual $n$-dependence in $\phi$, since
the fields $h_i$ are defined through $\bar q_N$, but this dependence is easily seen to be harmless 
thanks to the exponential decay of correlations inside
Dobrushin's uniqueness region, and to the fact that  $\bar q_N\to\bar q$.
Finally, defining $\tilde h_i=h+\beta \sqrt{\bar q}J_i$, the ergodic theorem implies that, for almost every $J$,
\begin{eqnarray*}
\lim_{n\to\infty}\sum_{k\in\vol} V^{n,k-1}(\xi_{n,k})
%%\stackrel{J-a.s.}{=}
=
\lim_{n\to\infty}\frac1n
\sum_{k\in\vol}\phi(\theta_{-k}\uh)=\EE\, \phi(\tilde{\uh})
\equiv \Gamma(\kappa,\beta,h),
\end{eqnarray*}
where we used the fact that the contribution of the spins $k\notin B_N$ vanishes for $|\vol|\to\infty$.
At this point, all the conditions necessary to 
 apply Theorem \ref{marty} are fulfilled and, recalling (\ref{riscalata}), one has (\ref{provided}),
which concludes the proof of Theorem \ref{th:2}.

\bigskip

{\bf Acknowledgments}

We would like to thank Giambattista Giacomin for many interesting discussions.
F.G. and F.L.T. are grateful to the {\sl Laboratoire de Probabilit\'es et Mod\`eles Al\'eatoires  P6 \& 7} 
for kind hospitality. F.L.T. was partially supported by
Swiss Science Foundation Contract No. 20-100536/1.

\small
%%%%%%%%%%%%%%%%%%

%%%%%%%%

\bigskip

\noindent
Francis Comets:
Universit{\'e} Paris 7 -- Denis Diderot, Math{\'e}matiques, case 7012, 2
place Jussieu, 
75251 Paris Cedex 05, France.
Email: comets@math.jussieu.fr\\
\texttt{http://www.proba.jussieu.fr/pageperso/comets/comets.html}\\

\noindent
Francesco Guerra:
Dipartimento di Fisica, Universit\`a di Roma ``La Sapienza'' and INFN, Sezione 
di Roma 1, P.le Aldo Moro 2, 00185 Roma, Italy. Email: francesco.guerra@roma1.infn.it\\

\noindent
Fabio Lucio Toninelli:
Laboratoire de Physique, UMR-CNRS 5672, ENS Lyon, 46 All\'ee d'Italie, 69364 Lyon Cedex 07, France.
Email: fltonine@ens-lyon.fr

\end{document}